\documentstyle[12pt]{article}
\newtheorem{theorem}{Theorem}

\newtheorem{corollary}{Corollary}

\newtheorem{conjecture}{Conjecture}

\newtheorem{claim}{Claim}
\newtheorem{problem}{Problem}

\newcommand{\text}[1]{\quad\mbox{#1}\quad}
\def\beq{\begin{equation}}\def\eeq{\end{equation}}
\def\beqn{\begin{eqnarray}}\def\eeqn{\end{eqnarray}}
\def\pont{\hspace{-6pt}{\bf.\ }}

\def\qed{\ifhmode\unskip\nobreak\fi\quad\ifmmode\Box\else$\Box$\fi}

\title{Rainbow matchings and partial transversals of Latin squares  \footnote { \noindent } }

\author{Andr\'as Gy\'arf\'as\\[-0.8ex]
\small Alfr\'ed R\'enyi Institute of Mathematics\\[-0.8ex]
\small Hungarian Academy of Sciences\\[-0.8ex]
\small Budapest, P.O. Box 127\\[-0.8ex]
\small Budapest, Hungary, H-1364 \small
\texttt{gyarfas.andras@renyi.mta.hu} \and G\'{a}bor N. S\'ark\"ozy
\thanks{Research supported in part by
the National Science Foundation under Grant No. DMS-0968699.}\\[-0.8ex]
\small Alfr\'ed R\'enyi Institute of Mathematics\\[-0.8ex]
\small Hungarian Academy of Sciences\\[-0.8ex]
\small Budapest, P.O. Box 127\\[-0.8ex]
\small Budapest, Hungary, H-1364
\small \texttt{sarkozy.gabor@renyi.mta.hu}
\\
\small Computer Science Department\\[-0.8ex]
\small Worcester Polytechnic Institute\\[-0.8ex]
\small Worcester, MA, USA 01609\\[-0.8ex]
\small \texttt{gsarkozy@cs.wpi.edu}}

\begin{document}
\maketitle

\begin{abstract}
In this paper we consider properly edge-colored graphs, i.e. two
edges with the same color cannot share an endpoint, so each color
class is a matching.  A matching is called \it rainbow \rm if its
edges have different colors. The minimum degree of a graph is
denoted by $\delta(G)$. We show that properly edge colored graphs
$G$ with $|V(G)|\ge 4\delta(G)-3$ have rainbow matchings of size
$\delta(G)$, this gives the best known estimate to a recent question
of Wang. Since one obviously needs at least $2\delta(G)$ vertices to
guarantee a rainbow matching of size $\delta(G)$, we investigate
what happens when $|V(G)|\ge 2\delta(G)$.

We show that any properly edge colored graph $G$ with $|V(G)|\ge
2\delta$ contains a rainbow matching of size at least $\delta -
2\delta(G)^{2/3}$. This result extends (with a weaker error term)
the well-known result that a factorization of the complete bipartite
graph $K_{n,n}$ has a rainbow matching of size $n-o(n)$, or
equivalently that every Latin square of order $n$ has a partial
transversal of size $n-o(n)$ (an asymptotic version of the Ryser -
Brualdi conjecture). In this direction we also show that every Latin
square of order $n$ has a {\em cycle-free partial transversal} of
size $n-o(n)$.
\end{abstract}

\section{Introduction - Rainbow matchings in proper colorings}

%The conjectures and results for factorizations have natural extensions for {\em proper edge colorings} i.e. edge colorings in which each color class has pairwise disjoint edges. L.D. Andersen \cite{A} proved that in any proper edge coloring of $K_n$ there is a rainbow linear forest (a forest with path components) with at least $n-\sqrt{2n}$ edges. The special case when the proper coloring is a factorization, it gives that every symmetric Latin square has a partial transversal with at least $n-\sqrt{2n}$ elements.

Recently Wang \cite{W} proposed to find the largest rainbow matching
in terms of the minimum degree in a properly colored graph. In fact,
\cite{W} raised the following problem.

%which is close to the Ryser - Brualdi conjecture in some sense.

\begin{problem}\label{wang} \pont Is it true that any properly colored graph $G$ contains a rainbow matching of size $\delta(G)$ provided that $|V(G)|$ is larger than a function of $\delta(G)$?
\end{problem}

Positive answers to Problem \ref{wang} were given in \cite{WZL},
\cite{DFMPW}, \cite{DFLMPW}, the current best bound is
${98\delta(G)\over 23}$ in \cite{DFLMPW}. In this paper we give a
better bound, namely $4\delta(G)-3$.

\begin{theorem}\pont \label{deltamatch} Any properly colored graph $G$ with at least $4\delta(G)-3$ vertices contains a rainbow matching of size $\delta(G)$.
\end{theorem}

Wang notices that the ``best'' function in his problem must be
greater than $2\delta(G)$ because certain Latin squares have no
transversals. For $\delta=2,3$ Theorem \ref{deltamatch} is best
possible, as shown by a properly $2$-colored $C_4$ for $\delta=2$
and by two vertex disjoint copies of a factorization of $K_4$ for
$\delta=3$. Our next result shows that if $n\geq 2\delta(G)$ then we
can find a rainbow matching almost as large as the desired
$\delta(G)$.

\begin{theorem}\pont\label{bigmatch1} Assume we have a proper coloring on a graph $G$ with $|V(G)|\ge 2\delta(G)$.
Then $G$ has a rainbow matching of size at least $\delta(G)- 2
(\delta(G))^{2/3}$.
\end{theorem}

Theorem \ref{bigmatch1} relates to partial transversals of Latin
squares. A {\em Latin square of order $n$} is an $n\times n$ array
$[a_{ij}]$ in which each symbol occurs exactly once in each row and
exactly once in each column. A {\em partial transversal} of a Latin
square is a set of distinct symbols, each from different rows and
columns. Latin squares can be also considered as {\em
$1$-factorizations of the complete bipartite graph $K_{n,n}$}, by
mapping rows and columns to vertex classes $R,C$ of $K_{n,n}$ and
considering the symbol $[a_{ij}]$ as the color of the edge $ij$ for
$i\in R, j\in C$. Then the edge sets with the same color form a
$1$-factorization of $K_{n,n}$ and partial transversals become
rainbow matchings. A well-known conjecture of Ryser \cite{RY} states
that for odd $n$ every $1$-factorization of $K_{n,n}$ has a rainbow
matching of size $n$. The companion conjecture, attributed to
Brualdi, is that for every $n$, every $1$-factorization of $K_{n,n}$
has a rainbow mathing of size at least $n-1$. These conjectures are
known to be true in an asymptotic sense, i.e. every
$1$-factorization of $K_{n,n}$ has a rainbow matching containing
$n-o(n)$ symbols. For the $o(n)$ term Woolbright \cite{WO} and
independently Brouwer at al. \cite{BVW} proved $\sqrt{n}$, Shor
\cite{SH} improved this to $5.518(\log{n})^2$ but it had an error
corrected in \cite{HS}. Theorem \ref{bigmatch1} extends these
results in two senses.  It allows proper colorings (instead of
factorizations) of arbitrary graphs (instead of complete bipartite
graphs).  On the other hand, the price we pay is that our error term
is weaker than the logarithmic one of Hatami and Shor \cite{HS}.

We also prove that Latin squares have a large partial transversal without short cycles.  A cycle of length $l$ in a Latin square $L$ means $a_{i_1j_1},a_{i_2j_2},\dots, a_{i_lj_l}$
such that $j_1=i_2,j_2=i_3\dots,j_l=i_1$ and all row (and column)
indices are different. For example, a cycle of length one is a
diagonal element of $L$, a cycle of length two is a pair of
symbols symmetric to the main diagonal, etc.

\begin{theorem}\pont\label{tetelcycle}
Assume that $L$ is a Latin square of order $n$ and $k\geq 2$ is a
positive integer.  Then $L$ has a partial transversal with at least
$n-6n^{\frac{k-1}{k}}$ elements containing no cycle of length $l$
for $l\le k$.
\end{theorem}

Applying Theorem \ref{tetelcycle}
with $$k=\left\lfloor \frac{\log n}{3\log\log n} \right\rfloor,$$
there is a partial transversal with at least
$n-6n^{\frac{k-1}{k}}$ elements that does not contain a cycle of length $l$ for
$l\leq k$. From each cycle (of length at least $k+1$) we remove an
arbitrary element of the transversal. In the resulting partial transversal we have at
least
$$n - \frac{n}{k+1} - 6n^{\frac{k-1}{k}}\geq
n - \frac{3n\log\log n}{\log n} - 6n^{\frac{k-1}{k}}\geq \left( 1 -
 \frac{4\log\log n}{\log n}\right)n $$ elements and we get the following.

\begin{corollary}\pont Any Latin square of order $n$ has a partial transversal $T$ of order $\left(1-\frac{4\log\log{n}}{\log{n}}\right) n$
such that $T$ has no cycles at all.
\end{corollary}

Notice that the error term in the corollary is much worse than in Theorem \ref{bigmatch1}. It is possible that the corollary holds in the following strong form (in the spirit of the Ryser - Brualdi conjecture).

\begin{conjecture}\label{cyclefree}\pont Any Latin square of order $n$ has a cycle-free partial transversal of order $n-2$.
\end{conjecture}
Conjecture \ref{cyclefree} would be best shown for $n=4$ by the Latin square $L$ with rows $1234,2143,3412,4321$. (One cannot select the symbol $1$ into a cycle-free partial transversal because it forms a loop, and only two of $\{2,3,4\}$ can be selected to avoid a $3$-cycle.)

\section{Proofs}

\subsection{Proof of Theorem \ref{deltamatch}}

Consider a properly colored graph $G$ with $|V(G)|\geq 4
\delta(G)-3$ and let $c(e)$ denote the color of edge $e$. We start
from a ``good'' configuration $H=M_1\cup M_2 \cup M_3 \cup_{i=1}^s
F_i$ defined as follows.
\begin{itemize}
\item For some integer $k\ge 0$  $M_1=\{e_i: i=1,2,\dots,k\}$ and $M_2=\{f_i: i=1,2,\dots,k\}$ form two vertex
disjoint rainbow matchings in $G$, $c(e_i)=c(f_i)$.

\item $M_3=\{g_i: i=k+1,\dots,\delta-1\}$ is a rainbow matching, vertex disjoint from
$M_1\cup M_2$ and $c(g_i)\ne c(e_j)$ for $1\le j \le k, k+1\le i\le
\delta-1$. Thus $M_1\cup M_3$ (likewise $M_2\cup M_3$) is a rainbow
matching of size $\delta -1$.

\item $F_1=\{h_i: i=k+1,\dots,t_1\}$ is a matching, vertex disjoint from $M_1\cup M_2$, and $h_i\cap M_3=\{v_i\}\in g_i$.  Moreover, $c(h_{k+1})\notin \{ c(e): e\in M_1\cup M_3\}$ and for $t_1\ge i>k+1$, $$c(h_i)\in \cup_{k+1\le j<i}c(g_j).$$ We call $F_1$ a chain. Note that $F_1$ is not necessarily rainbow, for example $c(h_i)=c(g_{k+1})$ for $k+1<i\le t_1$ satisfies the definition.

\item We allow several further disjoint chains $F_2,\dots,F_s$ where for $s\ge j\ge 2$, $F_j=\{h_i: i=t_{j-1}+1,\dots,t_j\}$ is a matching, vertex disjoint from $M_1\cup M_2\cup F_1\cup \dots\cup F_{j-1}$ and $h_i\cap M_3=\{v_i\}\in g_i$. Moreover, as before, $c(h_{t_{j-1}+1})\notin \{ c(e): e\in M_1\cup M_3\}$ and for
$t_j\ge i>t_{j-1}+1$, $$c(h_i)\in \cup_{t_{j-1}+1\le l<i}c(g_l)$$.
    \end{itemize}

One can easily see that a good configuration exists. Indeed, by
induction there is a rainbow matching $M$ with $\delta-1$ colors.
Let $v$ be a vertex not on $M$, and select an edge $vw$ of $G$ such
that $c(vw)\notin \{c(e): e\in M\}$. If $w$ is not on $M$ then $vw$
extends $M$ to a rainbow matching of size $\delta$ and the proof is
finished. Otherwise with $k=0,t_1=1$, $M_1=M_2=\emptyset,M_3=M,
F_1=\{vw\}=h_1$ we have a good configuration.

Select a good configuration $H$ with the largest possible $k$. Then select maximal chains $F_1,F_2,\dots,F_s$ to cover the maximum number of vertices of $M_3$ by $\cup_{i=1}^s F_i$.  If $k=\delta-1$, i.e. $M_3=F_1=\emptyset$  then select any vertex $v$ not in $H$ and an edge $vw$ such  that $c(vw)\notin \{ c(e): e\in M_1\}$. Since every color is repeated in $H$, we find a rainbow matching of size $\delta$. Thus we may assume $k<\delta-1$. Recall that  $v_i=g_i\cap h_i$ for $i=k+1,\dots t_s$.

Consider a vertex $v\notin V(H)$ and an edge $e=vw$ such that $c(vw)\notin \{c(f): f\in M_1\cup_{i>t_s} g_i$ and $w\ne v_i$ for $i=k+1,\dots t_s$. There is such an $e$ since we have precisely $\delta-1$ restrictions on the choice of $w$ and $\delta(v)\ge \delta$.

\noindent {\bf Case 1.} $w\in M_1\cup M_2$. If $j=c(vw)\notin
\{c(f): f\in M_1\cup M_3\}$ then if $w\in M_1$ (similarly if $w\in
M_2$) by adding the edge $vw$ to $M_2\cup M_3$, we find a rainbow
matching of size $\delta$. Otherwise from the choice of $w$,
$j=c(vw)=c(g_i)$ for some $t_s+1>i>k$.  We can now define a rainbow
matching of size $\delta$ as follows: for $1\le i \le k$ take either
$e_i$ or $f_i$ so that their union is disjoint from the edge $vw$.
This gives a matching with colors $1,\dots,k$ and color $j$. Remove
(the $j$-colored) $g_i$ from $M_3$ and add $h_i$ (from the chain
$F_l$ covering $v_i$). By definition of the chain, the color
$c(h_i)=c(g_{i_1})$ with $t_{l-1}+1\le i_1 <t_l$. Remove $g_{i_1}$
and add $h_{i_1}$ from $F_l$ and continue the procedure. Eventually
we end up by adding $h_{t_{l-1}+1}$  and the resulting matching is a
rainbow matching of size $\delta$.

\noindent {\bf Case 2.} $w\in \cup_{i=1}^s F_i$. For $c(vw)\in \{c(f): f\in M_1\cup M_3\}$ this contradicts the choice of $k$ since $vw$ can be added to the matching $M_1\cup M_2\cup M_3$ to get a new repeated color. For $c(vw)\notin \{c(f): f\in M_1\cup M_3\}$ we can add $vw$ to $M_1\cup M_3$ to get a rainbow matching of size $\delta$.

\noindent {\bf Case 3.} $w\in \cup_{i=k+1}^{t_s} g_i$, say $g_i=wv_i$ (by the choice of $c(vw)$, $w\ne v_i$). Since $v_i$ is in some chain $F_l$, we can add the edge $vw$, delete $g_i$, add $h_i\in F_l$, repeatedly until we end up by bringing in the first edge of the chain $F_l$ which has color $p\notin \{c(f): \{f\in M_1\cup M_3\}$. Thus we either get a new good configuration with $\delta-1$ colors in which the color $c(vw)$ is repeated or a matching with at least $\delta$ colors. The latter case finishes the proof and the former contradicts the choice of $k$.

\noindent {\bf Case 4.} $w\in \cup_{i>t_s} g_i$. This contradicts
the maximality of the chain cover, either because $c(vw)\notin
\{c(f): f\in M_1\cup M_3\}$ when we can start a new chain, or
$c(vw)\in \{c(f): f\in M_3\}$ when we can continue an existing
chain.

Since the good configurations involved have at most $4(\delta -1)$
vertices and one further vertex $w$ is required to get the rainbow
matching of size $\delta(G)$ from it,  $|V(G)|\ge 4\delta(G)-3$ is
indeed a sufficient condition. \qed

\subsection{Proof of Theorem \ref{bigmatch1}}

Let $M_1 =\{e_1,\dots,e_k\}$ be a maximum rainbow matching in a
properly colored graph $G$. Assume indirectly that $k<
\delta(G)-2(\delta(G))^{2/3}$. Set $\delta=\delta(G),
R=V(G)\setminus V(M_1)$ and let  $C_1$ be the set
 of ``new'' colors, i.e. colors not used on $M_1$. We have
\beq\label{R} |R|> 4 \delta^{2/3}. \eeq

Select an arbitrary $v\in R$. Since $deg(v)\geq \delta$ and $M_1$ is
maximum, at least $\delta-k> 2 \delta^{2/3}$ edges go back from $v$
to $M_1$ in colors $C_1$: \beq\label{v-C_1}deg_{C_1}(v,V(M_1))> 2
\delta^{2/3}.\eeq Indeed, otherwise we could increase the size of
our matching $M_1$. This implies in particular that $\delta - 2
\delta^{2/3} > \delta^{2/3}$, i.e. \beq\label{fok}\delta^{1/3} >
3,\eeq and that for the number of edges in $C_1$ between $M_1$ and
$R$ we have the following lower bound \beq\label{elek}
\left|E_{C_1}(R, V(M_1))\right|
> 2 \delta^{2/3} |R|.\eeq

In order to define the sets $M_2$ and $C_2$ in the next iteration we
do the following. We classify the edges $e_i$ in $M_1$ into two
classes. We put $e_i=x_iy_i$ into $M_1'$ if and only if

\beq\label{sum} deg_{C_1}(x_i,R)+deg_{C_1}(y_i,R) \geq 4
\delta^{1/3}(>12),\eeq using (\ref{fok}).

We define $M_2=M_1\setminus M_1'$ and $C_2=C_1\cup \{c(e_i)\; | \;
e_i\in M_1'\}$, where again $c(e_i)$ denotes the color of edge
$e_i$. We have the following two crucial claims.

\begin{claim}\pont\label{claim-egy}
$|M_1'|\geq \frac{\delta^{2/3}}{2}$, i.e. $|M_2|\leq |M_1| -
\frac{\delta^{2/3}}{2}$.
\end{claim}

Indeed, otherwise using (\ref{R}) we get
$$\left|E_{C_1}(R,
V(M_1))\right| \leq |M_1'|(2|R|) + |M_2|(4 \delta^{1/3}) <
\delta^{2/3}|R| + 4 \delta^{4/3} < 2 \delta^{2/3}|R|,$$ in
contradiction with (\ref{elek}).

\begin{claim}\pont\label{claim-ketto}
For every vertex $v\in R$ we have
$$deg_{C_2}(v,V(M_2)) > 2 \delta^{2/3}.$$
\end{claim}

For the proof of this claim observe first that if $e_i=x_iy_i\in
M_1'$, then all $C_1$-edges incident to this edge must be incident
to one of the endpoints (say $x_i$ always) since otherwise we could
increase $M_1$ (using (\ref{sum})). Denote by $X_1$ the set of these
$x_i$ endpoints in $M_1'$ and by $Y_1$ the set of other endpoints.
Thus there is no $C_1$-edge between $Y_1$ and $R$ and for every
$x_i\in X_1$ there are at least $4 \delta^{1/3}$ $C_1$-edges from
$x_i$ to $R$.

Consider an arbitrary $v\in R$ and an edge $vw$ with $c(vw)\in C_2$.
First note that $w\not\in R$. Indeed, otherwise if $c(vw)\in C_1$,
then we could clearly increase $M_1$ and if $c(vw)=c(e_i)$ for some
$e_i\in M_1'$, then we could increase $M_1$ again by exchanging
$e_i$ with $vw$ and adding a $C_1$-edge from $x_i$ to a free
neighbor in $R$ (using (\ref{sum}) again).

Thus $w\in V(M_1)$. Next we show that $w\not\in Y_1$. Assume
otherwise that $w=y_j$ for some $y_j\in Y_1$. If $c(vw)\in C_1$,
then again we could increase $M_1$ by exchanging $e_j$ with $vw$ and
adding another $C_1$-edge from $x_j$ to a free neighbor in $R$ such
that this edge has a different color from $c(vw)$ (using
(\ref{sum})). If $c(vw)=c(e_i)$ for some $e_i\in M_1'$, then we
could increase $M_1$ again by deleting $e_i$ and $e_j$, adding $vw$
and adding one $C_1$-edge from $x_i$, one $C_1$-edge from $x_j$ to
free neighbors in $R$ such that the two edges have different colors.

Thus if $w\in M_1'$, then $w\in X_1$ and this implies Claim
\ref{claim-ketto}, since by using (\ref{v-C_1}) we get
$$deg_{C_2}(v,M_2) \geq deg_{C_1}(v,M_1) + |M_1'| - |M_1'|>
2\delta^{2/3}.$$

Suppose now that $M_j$ and $C_j$ are already defined for a $j\geq 2$
such that the two claims are true for $j$, i.e. \beq\label{M_j}
|M_j|\leq |M_{j-1}|-\frac{\delta^{2/3}}{2},\eeq and
\beq\label{v-C_j} deg_{C_j}(v,V(M_j)) > 2 \delta^{2/3}.\eeq In order
to define $M_{j+1}$ and $C_{j+1}$ we put the edges $e_i=x_iy_i\in
M_j$ into $M_j'$ if and only if
$$deg_{C_j}(x_i,R)+deg_{C_j}(y_i,R) \geq 4 \delta^{1/3}.$$
We define $M_{j+1}=M_j\setminus M_j'$ and $C_{j+1}=C_j\cup
\{c(e_i)\; | \; e_i\in M_j'\}$.

Then we have to show that the two claims remain true for $j+1$. The
proof of Claim \ref{claim-egy} for $j+1$ is identical (replacing indices $1,2$ by $j,j+1$). The proof of
Claim \ref{claim-ketto} for $j+1$ is also similar but we will have
longer exchange sequences. First we show again that if
$e_i=x_iy_i\in M_j'$, then all $C_j$-edges to $R$ must be incident
to one of the endpoints (say $x_i$ always). Assume otherwise that we
have two $C_j$-edges of different colors $x_iv_1$ and $y_iv_2$,
where $v_1,v_2\in R$.

We ``trace back" both edges to a $C_1$-edge. If $c(x_iv_1)\in C_1$
(and similarly for $y_iv_2$), then we are done. Otherwise, by
definition, there exists a $j_1<j$ such that there exists an edge
$x_{i_1}y_{i_1}\in M_{j_1}'$ with $c(x_iv_1)=c(x_{i_1}y_{i_1})$. We
find a $C_{j_1}$-edge $x_{i_1}v_{i_1}$ such that $v_{i_1}$ is a free
neighbor of $x_{i_1}$ in $R$. If $c(x_{i_1}v_{i_1})\in C_1$, then we
are done. Otherwise, we trace this edge back further until we can
find an edge $x_{i_s}y_{i_s}\in M_{j_s}'$ for which there is a free
$C_1$-neighbor $v_{i_s}$ of $x_{i_s}$ in $R$. We proceed similarly
for $y_iv_2$ but we always select $C_{j_t}$-edges in unused colors
to free vertices in $R$. At this point we can define an increased
rainbow matching $M^*$ from $M$ by deleting the edges $x_iy_i$,
$x_{i_t}y_{i_t}$ for $t\in \{1,\ldots , s\}$ and adding the edges
$x_iv_1$, $x_{i_t}v_{i_t}$ for $t\in \{1. \ldots ,s\}$ and similarly
for $y_iv_2$. Note that the above procedure goes through if the
number of available neighbors in $R$ is at least $2(j+1)$. Since the
number of available neighbors is at least $4\delta^{1/3}$, the above
works for $j+1$ as long as $j+1\leq 2\delta^{1/3}$. Let us denote
again the set of these $x_i$ endpoints in $M_j'$ by $X_j$ and the
set of other endpoints by $Y_j$. Thus there is no $C_j$-edge between
$Y_j$ and $R$ and for every $x_i\in X_j$ there are at least $4
\delta^{1/3}$ $C_j$-edges from $x_i$ to $R$.

Consider again an arbitrary $v\in R$ and an edge $vw$ with $c(vw)\in
C_{j+1}$. First we show again that $w\not\in R$. Again as above we
trace $vw$ back to a $C_1$-edge. If $c(vw)\in C_1$, then we are
done. Otherwise, as above we find a sequence of edges
$x_{i_t}y_{i_t}\in M_{j_t}'$, $x_{i_t}v_{i_t}$ with $v_{i_t}\in R$,
$t\in \{1,\ldots ,s\}$ such that
$$c(vw)=c(x_{i_1}y_{i_1})\in C_{j+1},
c(x_{i_t}v_{i_t})=c(x_{i_{t+1}}y_{i_{t+1}})\in C_{j_t}, t\in
\{1,\ldots s-1\}$$ and
$$c(x_{i_s}v_{i_s})\in C_1.$$
Then again we can define an increased rainbow matching $M^*$ from
$M$ by deleting the edges $x_{i_t}y_{i_t}$ for $t\in \{1,\ldots ,
s\}$ and adding the edges $vw$, $x_{i_t}v_{i_t}$ for $t\in \{1.
\ldots ,s\}$. Again this works for $j+1$ when $j+1\leq 2
\delta^{1/3}$.

Thus $w\in V(M_j)$. Finally we show again that $w\not\in Y_j$.
Assume otherwise that $w=y_i$ for some $y_i\in Y_j$. Then again as
above we can trace back this $vw$ edge to a $C_1$-edge and thus we
could increase our matching. Thus if $w\in M_j'$, then $w\in X_1$
and this implies Claim \ref{claim-ketto} for $j+1$ assuming $j+1\leq
2 \delta^{1/3}$, since by using (\ref{v-C_j}) we get
$$deg_{C_{j+1}}(v,M_{j+1}) \geq deg_{C_j}(v,M_j) + |M_j'| - |M_j'|>
2\delta^{2/3}.$$

However, applying Claims \ref{claim-egy} and \ref{claim-ketto} with
$l=\lfloor 2 \delta^{1/3}\rfloor$ we get
$$deg_{C_l}(v,V(M_l)) > 2 \delta^{2/3},$$
while
$$|M_l| \leq |M_1|-(l-1)\frac{\delta^{2/3}}{2}< \delta - 2 \delta^{2/3} - (2\delta^{1/3}-2)\frac{\delta^{2/3}}{2}= - \delta^{2/3}<0,$$
a contradiction. \qed

\subsection{Proof of Theorem \ref{tetelcycle}}

To prove Theorem \ref{tetelcycle} we use another translation of a
Latin square: we associate the symbol $[a_{ij}]$ as a color to the
edge $ij$ of the complete directed graph $\overrightarrow{K_n}$,
where we have a loop $ii$ at each vertex $i$ and between each pair
$\{i,j\}$ of distinct vertices we have two oriented edges, $ij,ji$.
In this representation edges of the same color form an $1$-regular
digraph, i.e. the union of vertex disjoint directed cycles. Thus
Latin squares of order $n$ are equivalent to {\em $1$-factorizations
of $\overrightarrow{K_n}$.}

Theorem \ref{tetelcycle} is equivalent to

\begin{theorem}\pont\label{tetel3dirfact}
For every positive integer $k\geq 2$, in any factorization of
$\overrightarrow{K_n}$, there is a rainbow subgraph with maximum
indegree and outdegree one with at least $n-6n^{\frac{k-1}{k}}$
edges that does not contain a directed cycle $\overrightarrow{C_l}$
with $l\leq k$.
\end{theorem}

{\bf Proof of Theorem \ref{tetel3dirfact}.} Consider a one-factorization
of a $\overrightarrow{K_n}$: each color class is a $1$-regular
directed graph. Subgraphs of $1$-regular digraphs will be called {\em linear digraphs.} We start from a rainbow linear $G_1$ digraph on $n$ vertices with $t$ edges that does not contain a directed cycle $\overrightarrow{C_l}$ with $l\leq k$
such that $t$ is maximum.

We will show that
$$t\geq n-6n^{\frac{k-1}{k}}.$$
Thus $G_1$ is a collection of directed cycles with length greater
than $k$, directed paths, and isolated vertices. We consider
isolated vertices as degenerate paths where the beginning point and
ending point of the path are the same. Following the orientations,
we can go ``forward" from every vertex of $G_1$ that is not isolated
or endpoint of a path.

We will define two nested sequences of sets $A_1\subset A_2 \subset
\ldots$ and $B_1\subset B_2\subset \ldots$. At each step $i\ge 2$,
with a slight abuse of notation, we shall define $A_i\setminus
A_{i-1}, B_i\setminus B_{i-1}$ as the set of ``new" vertices in
$A_i,B_i$. Define $A_1$ as the set of beginning vertices of the
paths, and $B_1$ as the set of end vertices of the paths. We clearly
have
$$|A_1|=|B_1|=n-t.$$
Consider the edges with the ending point in $A_1$, having of the
$n-t$ new colors (not used in $G_1$), denote the set of these edges
by $E_1$. We will identify some edges in $E_1$ as {\em forbidden}
edges. For a beginning vertex $u\in A_1$ of a path $P$ in $G_1$ the
edge $vu$ is forbidden if $v$ is a vertex on the path $P$ at a
distance $l$ from $u$, where $2\leq l \leq k-1$. Indeed, these edges
may {\em potentially} create short rainbow cycles which are not
allowed. Thus altogether we have at most $(k-2)(n-t)$ forbidden
edges in $E_1$. This implies that there is a new color (denoted by
$c_1$) that contains at most $f_1=k-2$ forbidden edges. Consider
those edges in $E_1$ which have color $c_1$ and remove the at most
$f_1$ forbidden edges. Denote the resulting edge set by $E_1^{c_1}$.
We have the following claim.
\begin{claim}\pont\label{ki}
Assume that $vu\in E_1^{c_1}$ with $u\in A_1$. Then $v\not\in B_1$.
\end{claim}
Indeed, otherwise we would get a rainbow linear subgraph with
$t+1$ edges that does not contain a $C_l$ with $l\leq k$ (since the
forbidden edges were removed), a contradiction with the fact that
$t$ was maximum.

Now we are ready to define $A_2$ and $B_2$. Since we have a
factorization, every vertex is the ending point of an edge colored
with $c_1$ and thus \beq\label{E_1}|E_1^{c_1}|\geq n-t-f_1.\eeq
Define
$$B_2\setminus B_1 = \{ v \; | \; vu\in E_1^{c_1}, u\in A_1 \}.$$
This indeed makes the ``new" vertices of $B_2$ disjoint from $B_1$
by Claim \ref{ki} and by (\ref{E_1}) we have $|B_2\setminus B_1|\geq
n-t-f_1$. To get $A_2\setminus A_1$ we shift the vertices in
$B_2\setminus B_1$ forward by one on their paths or cycles in $G_1$.
We can do this shifting since the vertices in $B_2\setminus B_1$ are
never isolated or endpoints of the paths (those vertices are in
$B_1$). Furthermore, clearly this set is indeed disjoint from $A_1$
since we are shifting away from the beginning vertices of the paths.
Finally, define $G_2$ as the union of $G_1$ and the edge set
$E_1^{c_1}$. Assume that there is a rainbow $C_l$ in $G_2$ with
$l\leq k$. Then this $C_l$ contains exactly one edge colored with
$c_1$ (since $G_1$ does not contain a rainbow $C_l$ with $l\leq k$),
but then this edge is forbidden and was removed, a contradiction.

At this point we have the following three properties for $i=2$ (with
$f_1=k-2$):

\begin{enumerate}
\item $A_{i-1}\subset A_i$ and $B_{i-1}\subset B_i$,
\item $n-t\geq |A_i\setminus A_{i-1}|=|B_i\setminus B_{i-1}|\geq n-t-f_{i-1}$,
\item $G_i$ does not contain a rainbow $C_l$ with $l\leq k$,
\item  For every $u\in A_i$ there is a linear rainbow subdigraph of $G_i$ with $t$ edges
such that the set of endpoints of the paths is exactly $B_1$ and one
of the paths begins at $u$.
\end{enumerate}

We continue in this fashion, maintaining these properties with a
suitable $f_{i-1}$. Assume that $A_1, A_2, \ldots , A_{i-1}$ and
$B_1, B_2, \ldots , B_{i-1}$ are already defined for some $i\geq 3$.
%We will define $A_i$, $B_i$ and $G_i$ such that
%\begin{enumerate}
%\item $A_{i-1}\subset A_i$ and $B_{i-1}\subset B_i$,
%\item $n-t\geq |A_i\setminus A_{i-1}|=|B_i\setminus B_{i-1}|\geq n-t-\frac{2ki^{k-1}(n-t)}{n-t-(i-2)}$,
%\item $G_i$ does not contain a rainbow $C_l$ with $l\leq k$.
%\end{enumerate}
Consider the edges with the ending point in $A_{i-1}$ in one of the
$n-t-(i-2)$ new colors (not used in $G_{i-1}$), denote the set of
these edges by $E_{i-1}$. Again we will identify some edges in
$E_{i-1}$ as forbidden edges. For a vertex $u\in A_{i-1}$ the edge
$vu\in E_{i-1}$ is forbidden if there is a rainbow path of length at
most $k-1$ from $u$ to $v$ in $G_{i-1}$, where the last edge is from
$G_1$. Indeed, these edges may potentially create short rainbow
cycles which are not allowed. For a fixed $u\in A_{i-1}$ for the
number of these rainbow paths of length at most $k-1$ (and thus for
the number of $vu\in E_{i-1}$ forbidden edges) a crude upper bound
is $(k-2) i^{k-2}$. Indeed, for each of the edges before the last
one we have at most $i$ possibilities (one from $G_1$ and one for
each of the $i-2$ added colors) and for the last edge we have only
one possibility as it must be in $G_1$.

Thus altogether we have at most $ki^{k-1}(n-t)$ forbidden edges in
$E_{i-1}$. This implies that there is a new color (denoted by
$c_{i-1}$) that contains at most
\begin{equation}\label{star}
f_{i-1}=\frac{ki^{k-1}(n-t)}{n-t-(i-2)}
\end{equation} forbidden
edges. Consider those edges in $E_{i-1}$ which have color $c_{i-1}$
and remove these forbidden edges. Denote the resulting edge set by
$E_{i-1}^{c_{i-1}}$. We have the following claim.
\begin{claim}\pont\label{kii}
Assume that $vu\in E_{i-1}^{c_{i-1}}$ with $u\in A_{i-1}$. Then
$v\not\in B_1$.
\end{claim}
Otherwise from property 4 the edge $vu$ would join two paths or create a cycle and we would get
a rainbow linear subgraph with $t+1$ edges, a contradiction with the fact that $t$ was
maximum. We would have no rainbow $C_l$ with $l\leq k$ either since
otherwise this $C_l$ must contain exactly one edge colored with
$c_{i-1}$, namely $vu$ (since $G_{i-1}$ does not contain a rainbow
$C_l$ with $l\leq k$), but then this edge is forbidden and was
removed, a contradiction.

Now we are ready to define $A_i$ and $B_i$. Since we have a
factorization, every vertex is the ending point of an edge colored
with $c_{i-1}$. Define
$$B_i\setminus B_{i-1} = \{ v \; | \; vu\in E_{i-1}^{c_{i-1}}, u\in A_{i-1}, v\not\in B_{i-1} \}.$$

This indeed by definition makes the ``new" vertices of $B_i$
disjoint from $B_{i-1}$ and by Claim \ref{kii} we have property 2
for $|B_i\setminus B_{i-1}|$, since
$$|B_i\setminus B_{i-1}| \geq |A_{i-1}|-f_{i-1} - |B_{i-1}\setminus
B_1|= |A_{i-1}|-f_{i-1} - |A_{i-1}\setminus A_1| =$$ $$= |A_1|-
f_{i-1} = n-t-f_{i-1}.$$ To get $A_i\setminus A_{i-1}$ we shift the
vertices in $B_i\setminus B_{i-1}$  by one step forward on their
paths or cycles in $G_1$. We can do this shifting since the vertices
in $B_i\setminus B_{i-1}$ are never isolated or endpoints of the
paths (those vertices are in $B_1$). Furthermore, this set is indeed
disjoint from $A_{i-1}$ since in $G_1$ the in-degree of every vertex
is at most one.

Finally, define $G_i$ as the union of $G_{i-1}$ and those edges in
$E_{i-1}^{c_{i-1}}$ that start in vertices in $B_i\setminus
B_{i-1}$. Notice that property 4 is maintained from the definition of the
``new" vertices of $A_i\setminus A_{i-1},B_i\setminus B_{i-1}$.

Property 3 above is also true for $G_i$. Indeed,
otherwise assume indirectly that there is a rainbow $C_l$ in $G_i$
with $l\leq k$. Then this $C_l$ contains exactly one edge colored
with $c_{i-1}$ (since $G_{i-1}$ does not contain a rainbow $C_l$
with $l\leq k$), but then this edge is forbidden and was removed, a
contradiction. Note again that from the construction indeed the last
edge is from $G_1$ on the rainbow path of length at most $k-1$
connecting the two endpoints of this edge.

Next we claim that \beq\label{expand}|A_i\setminus
A_{i-1}|=|B_i\setminus B_{i-1}|\geq \frac{n-t}{2} \; \; \mbox{for}\;
\;  i\leq \left(\frac{n-t}{4k}\right)^{\frac{1}{k-1}}.\eeq In fact,
for these $i$'s from property 2 we have $$|A_i\setminus
A_{i-1}|=|B_i\setminus B_{i-1}|\geq
n-t-\frac{ki^{k-1}(n-t)}{\frac{n-t}{2}} = n - t - 2ki^{k-1} \geq
\frac{n-t}{2}.$$ Thus we must have
$$\frac{n-t}{2} \left(\frac{n-t}{4k}\right)^{\frac{1}{k-1}} \leq
n,$$ and therefore using $k\geq 2$
$$n-t \leq 2 (4k)^{\frac{1}{k}} n^{\frac{k-1}{k}}\leq 6 n^{\frac{k-1}{k}}.$$
From this we get that
$$t \geq n - 6 n^{\frac{k-1}{k}},$$ as desired. $\Box$

\end{document}